%% file: 010_main.tex
\documentclass{amsart}

\input{000_setup}

\title{Average Elliptic Billiard\\ Invariants with Spatial Integrals}
\author{Jair Koiller}
\author{Dan Reznik}
\author{Ronaldo Garcia}

\date{January 2020}

\begin{document}

\maketitle

\input{005_abstract}

\section{Introduction}
\label{sec:introduction}
\input{020_introduction}

\section{Preliminaries}
\label{sec:preliminaries}
\input{025_preliminaries}

\section{Average Sidelength}
\label{sec:perimeter}
\input{030_perimeter}

\section{Average Cosine}
\label{sec:cosines}
\input{040_cosines}

\section{Geometric Mean of Outer Cosines}
\label{sec:outer-cosines}
\input{060_outer_cosines}


\section{Questions}
\label{sec:questions}
\input{080_questions}

\section*{Acknowledgments}
\input{110_acknowledgments}

\bibliographystyle{maa}
\bibliography{references,authors_rgk_v3,jair_spatial_v1}

\end{document}

%% file: 000_setup.tex
\usepackage{graphics}
\usepackage{thmtools}
\usepackage{wasysym}
\usepackage[T1]{fontenc}    

\usepackage{amsthm}
\usepackage{amsbsy,amsmath,amssymb,amscd,amsfonts}
\usepackage[pagebackref=true]{hyperref}
\usepackage[nameinlink,capitalize,noabbrev]{cleveref}

\usepackage{graphicx,float,latexsym,color}
\usepackage[font={scriptsize,it}]{caption}
\usepackage{subcaption}

\usepackage{makecell}

\usepackage[dvipsnames]{xcolor}

\newtheorem*{theorem*}{Theorem}

\theoremstyle{remark}

\theoremstyle{definition}

\hypersetup{
    pdftoolbar=true,        
    pdfmenubar=true,        
    pdffitwindow=false,     
    pdfstartview={FitH},    
    colorlinks=true,       
    linkcolor=OliveGreen,          
    citecolor=blue,        
    filecolor=black,      
    urlcolor=red           
}

\usepackage{lineno}

\usepackage{comment}
\usepackage{microtype}
\usepackage{footnote}

\newcommand{\E}{\mathcal{E}}

%% file: 005_abstract.tex
\begin{abstract}
We compare invariants of N-periodic trajectories in the elliptic billiard, classic and new, to their aperiodic counterparts via a spatial integrals evaluated over the boundary of the elliptic billiard. The integrand is weighed by a universal measure equal to the density of rays hitting a given boundary point. We find that aperiodic averages are smooth and monotonic on caustic eccentricity, and perfectly match N-periodic average invariants at the discrete caustic parameters which admit a given N-periodic family. 
\end{abstract}

%% file: 020_introduction.tex
\begin{figure}
 \centering
 \includegraphics[width=\textwidth]{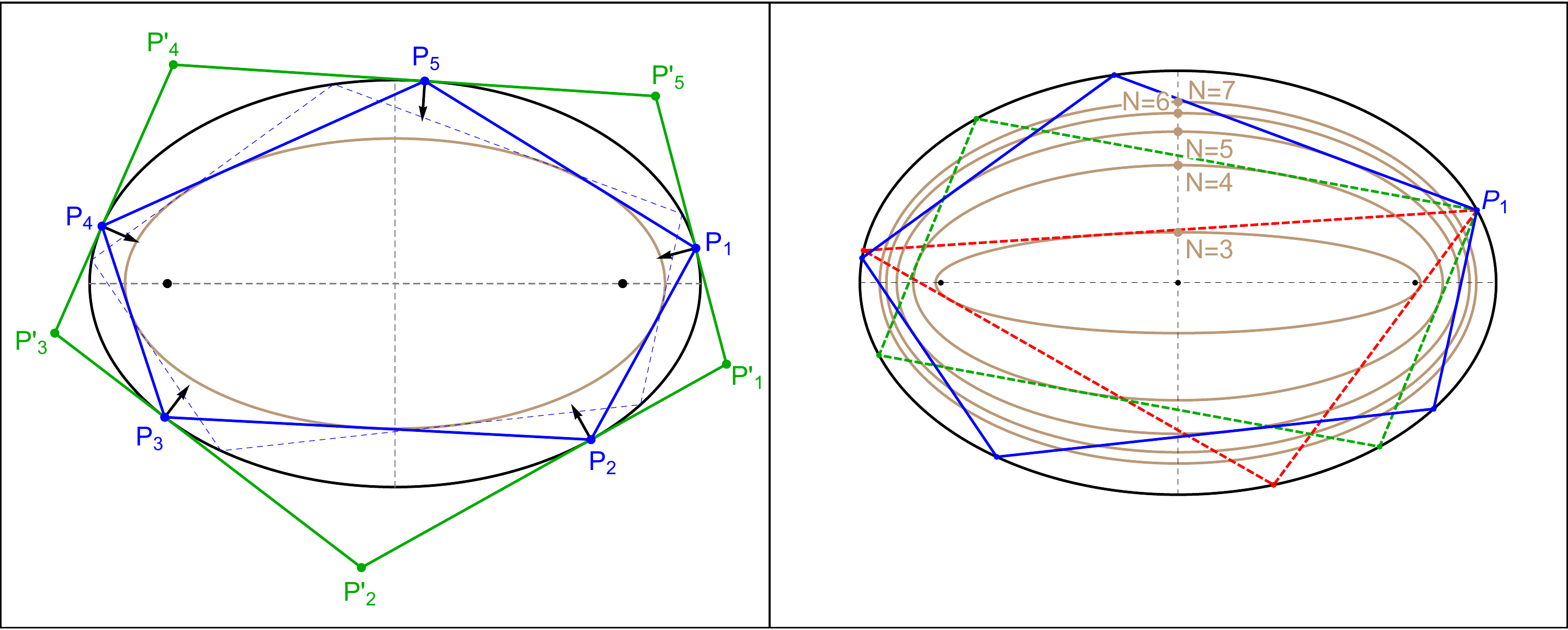}
 \caption{\textbf{Left:} A 5-periodic trajectory (blue) in the elliptic billiard (black), whose vertices are bisected by the ellipse normals (black arrows). The Poncelet family remains tangent to a confocal elliptic caustic (brown). A second, same perimeter 5-periodic is also shown (dashed blue). The outer polygon (green) has sides tangent to the elliptic billiard at the vertices of the N-periodic. \textbf{Right:} discrete confocal caustics (brown) associated with N-periodics, N=3,...,7. Show are sample 3-, 4-, and 5-periodics (dashed red, dashed green, solid blue, respectively) sharing one common vertex $P_1$.}
 \label{fig:definitions}
\end{figure}

The two classic invariants of Poncelet N-periodics in the elliptic billiard are perimeter $L$ and quantity known as Joachimsthal's constant $J$; see \cref{fig:definitions}. The former implies a billiard trajectory is an extremum of the perimeter function while the latter is equivalent to stating that all trajectory segments are tangent to a confocal caustic \cite{sergei91}.

Experiments have unearthed a few additional ``dependent'' invariants\footnote{A billiard N-periodic is fully specified by $L,J$, so any ``new'' invariants are ultimately dependent on them.} including (i) the sum of cosines, (ii) the product of outer polygon cosines, (iii) certain ratios of areas, etc. \cite{reznik2020-intelligencer}. These have been subsequently proved \cite{akopyan2020-invariants, bialy2020-invariants,caliz2020-area-product}. More recently, the list of conjectured invariants has grown to many dozen \cite{reznik2021-invariants}.

\begin{figure}
 \centering
 \includegraphics[width=\textwidth]{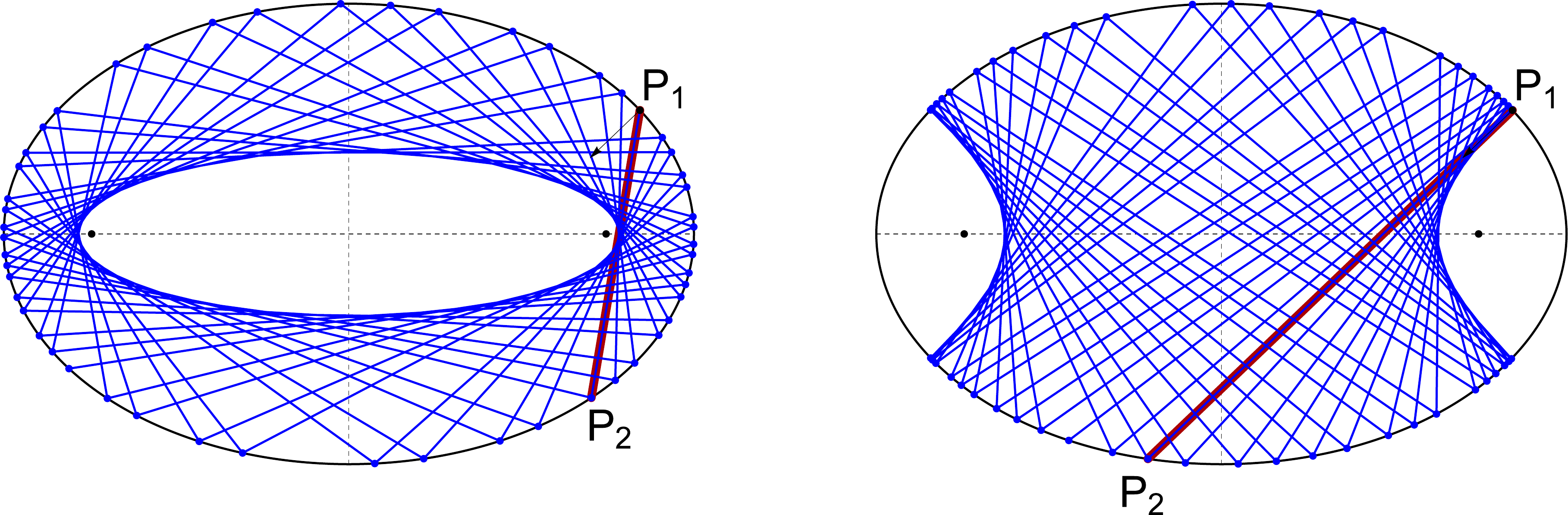}
 \caption{Two regimes of aperiodic, space-filling trajectories in an elliptic billiard, reproduced from \cite{reznik2020-intelligencer}. \textbf{Left:} initial ray $P_1 P_2$ does not pass between the foci, confocal caustic is ellipse. \textbf{Left:} initial ray passes between the foci, confocal caustic is hyperbola.}
 \label{fig:aperiodic}
\end{figure}

With a small perturbation of the caustic, an N-periodic trajectory becomes aperiodic (space-filling); see \cref{fig:aperiodic}. A key question we explore is: given a discrete invariant computed for an N-periodic, what is its analogue in the space-filling case? In the latter case, the sum or product of a given quantity can diverge. Fortunately, in both cases we can compare their finite {\em averages}. 

\subsection*{Main Result}

Our contribution is to accurately and efficiently estimate aperiodic averages using a {\em spatial integral} evaluated over the caustic's boundary. We weigh the integrand by a the elliptic billiard universal measure \cite[Section 51]{arnold-book78}, \cite{poritsky50,lazutkin73,jovanovic11-hamilton} which yields the aperiodic density of rays hitting a particular point on the billiard boundary.

Referring to \cref{tab:invs}, we will examine one classic (perimeter) and two ``new'' invariants (sum of cosines and product of exterior cosines). We will compare their averages (average chord length, average cosine, and geometric mean of exterior cosines) within the continuum of aperiodic trajectories.

\begin{table}
\begin{tabular}{|c|c|c|c|}
\hline
invariant & formula & from & average \\
\hline
L & elliptic functions & classic & $L/N$ \\
$\sum{\cos\theta_i}$ & $L J - N$ & \cite{reznik2020-intelligencer} & $(LJ)/N-1$\\
$\prod{\cos\theta_i'}$ & ? & \cite{reznik2020-intelligencer} & geometric mean \\
\hline
\end{tabular}
\caption{Three N-periodic invariants whose averages are compared to those displayed by their aperiodic counterparts. A ``?'' means no closed expression has yet been derived.}
\label{tab:invs}
\end{table}

\subsection*{Article Structure}
In \cref{sec:preliminaries} we review preliminary concepts and definitions. In \cref{sec:perimeter,sec:cosines,sec:outer-cosines} we calculate, via spatial integrals, the (i) average sidelength (i.e., average perimeter), (ii) average cosine, and (iii) geometric mean of outer cosines. We then compare them with the values predicted by either closed form or numeric computation of the original quantities in the N-periodic case, showing that they lie at consistent locations within the continuum of confocal caustics. Unanswered questions and/or future work appear in \cref{sec:questions}.

%% file: 025_preliminaries.tex
Let $(\E,\E_c)$ denote the outer and inner ellipses in the confocal pair given by:

\[\E:\;\frac{x^2}{a^2} +  \frac{y^2}{b^2} = 1,\;\;\;
    \E_c:\;\frac{x^2}{a_c^2} +  \frac{y^2}{b_c^2} = 1
\]

Let $c^2=a^2-b^2=a_c^2-b_c^2$. Let $\mathcal{A}=\mbox{diag}[1/a^2,1/b^2]$. Joachimsthal's constant at a point $\delta$ on $\E$ is given by $\left<\mathcal{A}.P,v\right>$, where $v$ is the unit incoming vector \cite{sergei91}. $J$ is also given by:

\begin{equation*}
    J = \frac{\sqrt{\lambda}}{{a}{b}}
\end{equation*}

\noindent where $\lambda=a^2 - a_c^2=b^2-b_c^2$.

\begin{figure}
    \centering
    \includegraphics[width=.7\textwidth]{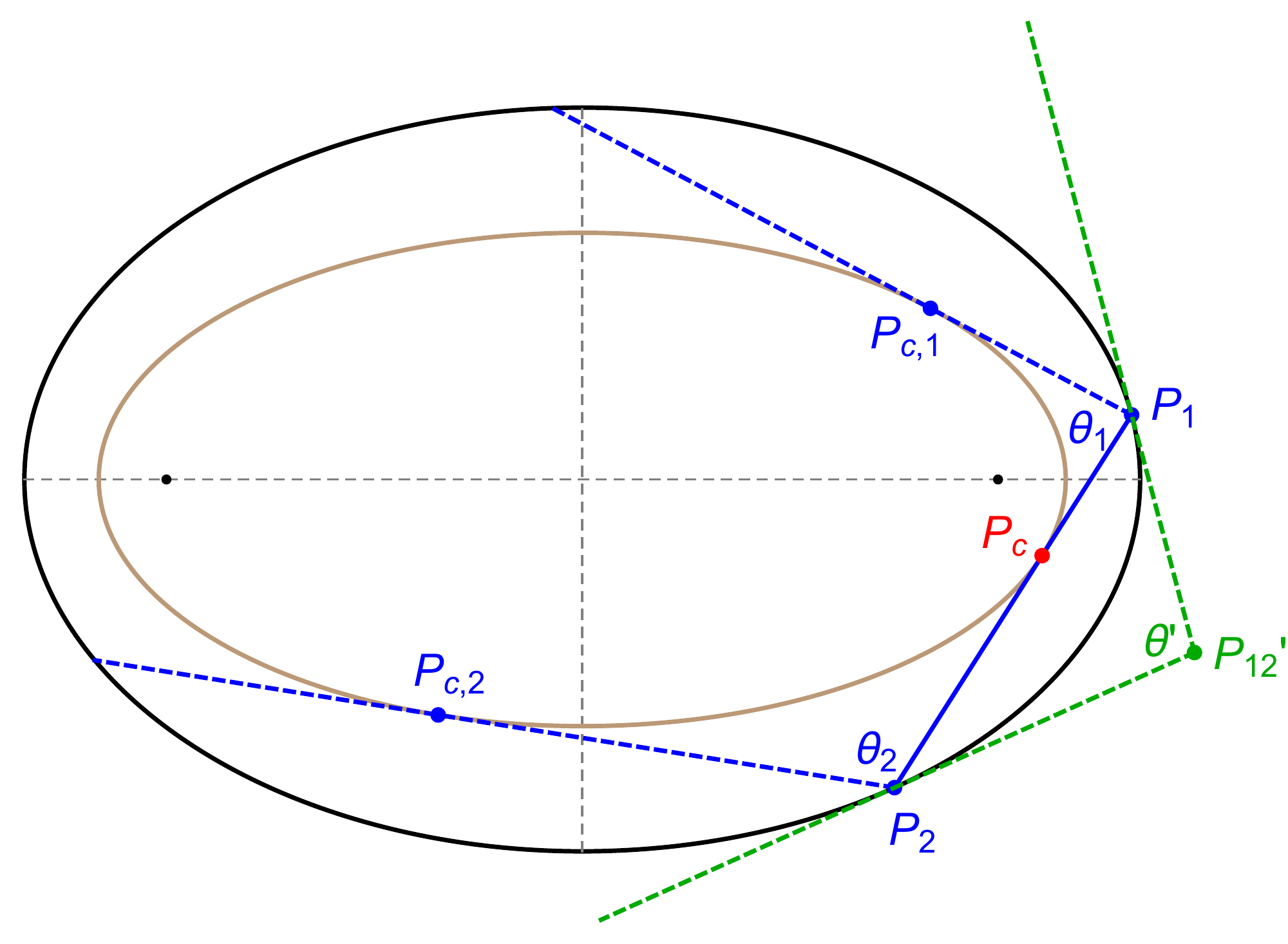}
    \caption{The tangent to a chosen point $P_c$ on the caustic intersects the elliptic billiard at $P_1$ and $P_2$. Let $\theta_1$ (resp. $\theta_2$) be the angle between segment $P_1 P_2$ and the next $P_{c,1}$ (resp. previous $P_{c,2}$) tangent to the caustic (dashed blue). The outer angle $\theta'$ associated with $P_c$ is measured at the intersection $P_{12}'$ of the tangents to the billiard (dashed green) at $P_1$ and $P_2$.}
    \label{fig:p12}
\end{figure}

Referring to Fig.~\ref{fig:p12}, let $P_c$ be a point on the caustic, and $P_1,P_2$ be the intersections of the tangent through $P_c$ with the outer ellipse. These are given by:

\begin{align}
P_1&=(x_1,y_1)=\frac{1}{\psi}\left[a_c^2a(ab_c^4x_c-\zeta by_c),b_c^2b(ba_c^4y_c+\zeta ax_c)\right] \label{eqn:p1} \\
P_2&=(x_2,y_2=\frac{1}{\psi}\left[a_c^2a(ab_c^4x_c+\zeta by_c),b_c^2b(ba_c^4y_c-\zeta ax_c)\right] \label{eqn:p2}
\\
\zeta&=\sqrt{a^2-a_c^2}\sqrt{ b_c^4   x_c^2 +   a_c^4  y_c^2  } \nonumber\\
\psi=&a^2b_c^4 x_c^2+b^2a_c^4y_c^2\nonumber
\end{align}

Given a confocal caustic (say parametrized by its minor axis $b_c$) let $P_i$ be the vertices of an associated aperiodic trajectory. $i=1,...,\infty$. The asymptotic average $\overline{g}$ of some vertex-evaluated quantity $g(P_i)$ is given by:

 \begin{equation}
 \overline{g} =  \lim_{N \to \infty}\;\frac{1}{N}\;\sum_{i=1}^k g(P_i)
 \label{eqn:asympt-101}
 \end{equation}

One can evaluate $g(s)$, for example at either intersection $P_1$ or $P_2$ in Figure~\ref{fig:p12}.

The billiard map is an involution of the pair $(P,P_c)$ of a point on $\E$ and $\E_c$ respectively to new points There is change of variables $s \to x$ which linearizes the billiard map, $x \to  x + \tau$ \cite[Chapter 13]{flatto2009}, \cite{kolod1985,glutsyuk19-poritsky,zhang17}.



Let $\rho$ be the density of an invariant measure normalized such that $\oint \rho(x)ds=1$. This can be regarded as the density of rays associated with $x$. The following universal measure has been derived for the elliptic billiard, independent of $\tau$:  

\begin{equation}
dx = \kappa_c^{2/3} ds
\end{equation} 

The above allows us to replace \eqref{eqn:asympt-101} with the following spatial integral:

\begin{equation}
\overline{g} = \frac{1}{\oint\,\kappa_c^{2/3}ds}\,\,\oint\,g(s)\ \kappa_c^{2/3} ds  .
\label{eqn:spatial-general}
\end{equation}

\subsection*{Auxiliary expressions}
 
Arc length and curvature along the caustic ellipse are given by:

\begin{align*}
ds = & \left( a_c^2 \, \sin^2 u + b_c^2 \,  \cos^2 u       \right)^{1/2} du\\
\kappa_c = &  (a_c b_c)^{2/3} \, \left( a_c^2 \, \sin^2 u + b_c^2 \,  \cos^2{u} \right)^{-3/2}
\end{align*}

\noindent so that:

\begin{align*}
\rho = dx = &  \kappa_c^{2/3} ds \\
= & (a_c b_c)^{2/3} \, \left( a_c^2 \, \sin^2 u + b_c^2 \,  \cos^2 u       \right)^{-1/2} \\
 = & \frac{(a_c b_c)^{2/3}}{\sqrt{a_c^2 y^2/b_c^2 + b_c^2 x^2/a_c^2}}=\frac{(a^2-\lambda)^{\frac{1}{3}}(b^2-\lambda)^{\frac{1}{3}}}{\sqrt{a^2-\lambda-(a^2-b^2)\cos^2u }}
\end{align*}

Below we will be also expressing certain average quantities in terms of the following Jacobi elliptic functions of the first and third kind, respectively \cite[Introduction]{grad1980}:

 \begin{align*}
  K(m)=&\int_0^{\frac{\pi}{2}} \frac{d\alpha}{\sqrt{1-m^2\sin^2\alpha}}\\
 \Pi(n,m)=&\int_0^{\frac{\pi}{2}}\frac{d\alpha}{(1-n^2\sin^2\alpha)\sqrt{1-m^2\sin^2\alpha}  },\;\; 
 \end{align*}

%% file: 030_perimeter.tex
We constuct a spatial integral to compute $\overline{L}$, the average  sidelength in an aperiodic trajectory and compare it with $L/N$ for an N-periodic.

The distance between two consecutive points $P_1$ and $P_2$ of a billiard orbit parametrized by the point $P_c=[x_c,y_c]$ in the confocal caustic is given by:

\[  l_{12}=\,{\frac {2ab\sqrt {    {a}^{2} b_c^{4}x_c^2 
   +a_c^{4}b^2y_c^2 -a_c^{4}b_c^{4}  }\sqrt {a_c
 ^{4}y_c^2 +b_c^{4}x_c^2 }}{{a}^{2}b_c^{4}x_c^2 +a_c^{4}{b}^{2}y_c^2 }}\]
Therefore, 
$P_c=[\sqrt{a^2-\lambda}\cos u, \sqrt{b^2-\lambda}\sin u]$ leads to
\[l_{12}(u)= 2ab\sqrt {\lambda}\,\frac {    c^2 \cos  ^{2}u  -   ({a}^{2}-\lambda)
  }{  \lambda c^2 \cos^{2}u  
 -({a}^{2}-\lambda){b}^{2} }
\]

\begin{align*}
 l(u) &=l_{12}\kappa_c^{2/3}ds=\,{\frac {c_1\sqrt {{a}^{2}-\lambda-{c}^{2}  \cos^2u  }}{{b}^{2} \left( {a}^{2}-
\lambda \right)-\lambda{c}^{
2}    \cos^2u    }}\\
c_1&=2ab\sqrt {\lambda} \sqrt [3]{{a}^{2}-\lambda}\sqrt [3]{{b}^{2}-\lambda}
\end{align*}
\begin{equation}
\label{eqn:Lmed}
\overline{L} = \frac{1}{\oint \kappa_c^{2/3}  ds}  \int_0^{2\pi} \, l(u)du
\end{equation}

\begin{figure}
    \centering
    \includegraphics[width=\textwidth]{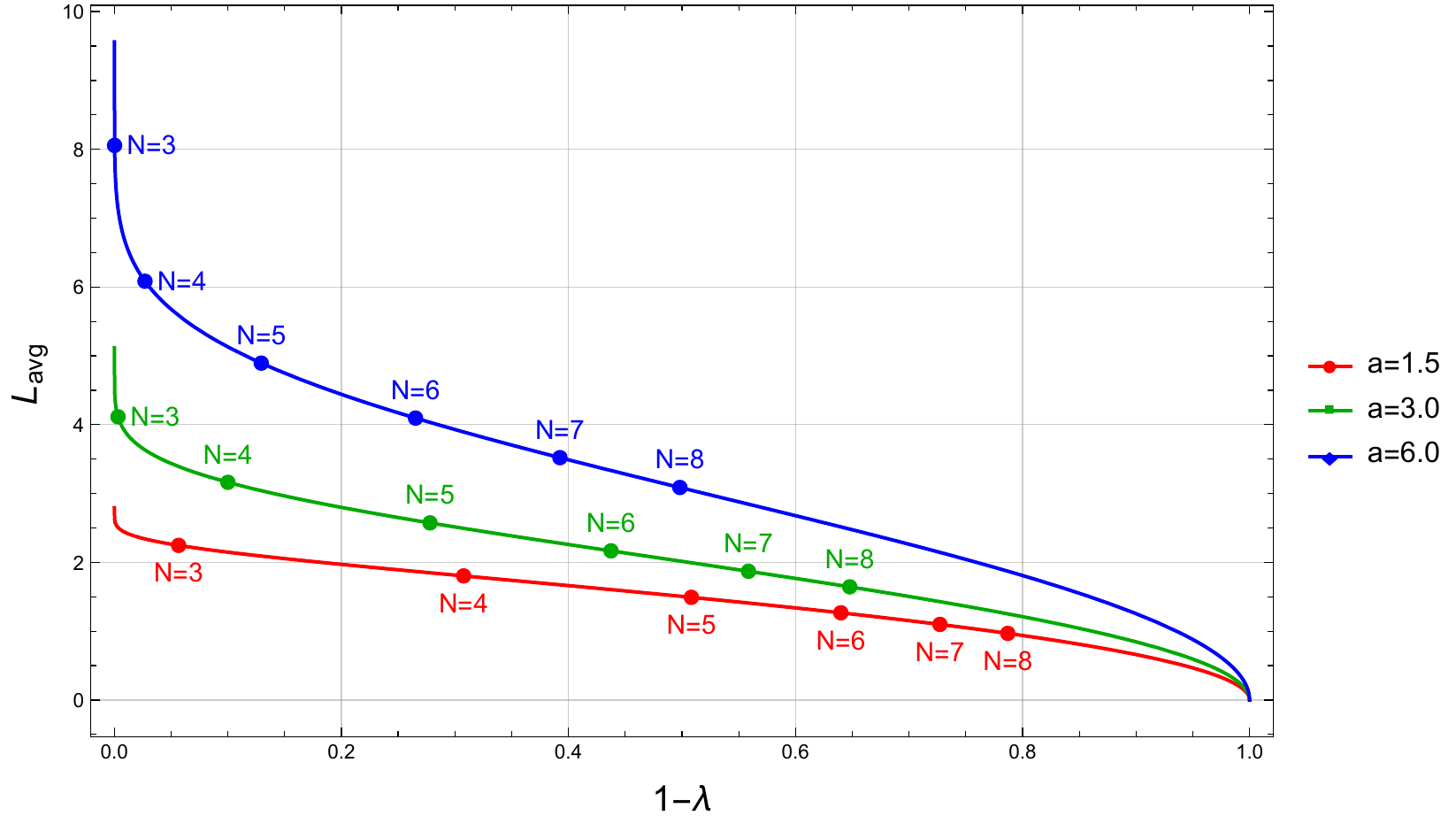}
    \caption{The value of the average  sidelengths vs $1-\lambda$, $b=1$, and three values of $a$. The dots show agreement of the value with $L/N$ for various non-intersecting N-periodics. When $1-\lambda$ is zero, the average perimeter tends to $2a$.}
    \label{fig:lmed-lambda}
\end{figure}

In terms of the elliptic integrals $K$ and $\Pi$ we have that:

\begin{align*}
\int {\kappa_c}^{\frac{2}{3}}ds&=\int_0^{2\pi} \frac{(a^2-\lambda)^{\frac{1}{3}}(b^2-\lambda)^{\frac{1}{3}}}{\sqrt{s_3} \sqrt{1-s_3\cos^2u}}du\\
&= \frac{4(a^2-\lambda)^{\frac{1}{3}}(b^2-\lambda)^{\frac{1}{3}}}{\sqrt{s_3}} K(\sqrt{s_3}) \\
\int_0^{2\pi} l(u)du&=\frac{c_1}{b^2} {\frac { \left(   2\,s_5 -2s_3\right) \Pi
 \left( s_5,\sqrt {s_3} \right) +2\,s_3 K\left( \sqrt {s_3} \right) }{s_5}}\\
s_3&=\frac{c^2}{a^2-\lambda}, \;\;s_5=\frac{\lambda s_3}{b^2}  
\end{align*}
Therefore,

\[\overline{L}= {\frac {2a}{\sqrt {\lambda}\,K\left( \sqrt {s_3}\right) } \left(  \left( -{b}^{2}+\lambda \right) \Pi \left( {\frac {\lambda\,s_3}{{b}^{2}}},\sqrt {s_3}
 \right) +K \left( \sqrt {s_3} \right) {b}^{2}
 \right) }
\]
Numerical results are shown in \cref{fig:lmed-lambda} for three different billiard aspect ratios. Notice points on each curve report the average perimeters $L/N$ obtained with N-periodics at the required caustic parameters $\lambda$.

%% file: 040_cosines.tex
We evaluate the average cosine $\overline{C}$ for aperiodics with a spatial integral. 

Using the Joachimstall invariant we 
obtain:

 \begin{equation} \aligned
 \cos\theta_{1}  &= \frac{J^2a^4b^4}{ 2(a^4y_{1}^2+b^4x_{1}^2)}-1=\frac{\lambda a^2 b^2}{2(a^4 y_{1}^2+b^4 x_{1}^2)}-1\\
 &=\frac{\lambda}{2d_1d_2}-1,\;\; d_1 =|P_1-F_1|,\;\; d_2=|P_1-F_2|
 \endaligned
 \end{equation}



 
\noindent Let $\cos\theta(u)=(\cos\theta_{1}(u)+\cos\theta_{2}(u))/2$. Let $a_c=\sqrt{a^2-\lambda}$, $b_c=\sqrt{b^2-\lambda}$
and $(x_c,y_c)=(a_c\cos u,b_c\sin u)$.

 Using \eqref{eqn:p1}  and \eqref{eqn:p2} 
it follows that:

\begin{align*}
\cos\theta(u) &=  \frac{r_1+r_2\cos^2u}{r_3+r_4\cos^2u}= \frac{r_1}{r_3}\frac{1-s_1\cos^2 u}{1-s_2\cos^2u}\\
r_1&=(a^2-\lambda)(a^2b^2-c^2\lambda )(2a^2b^2-2a^2\lambda+b^2\lambda)\\
r_2&=-c^2(a^2 b^2-(a^2+b^2)\lambda)(2a^2b^2-2a^2\lambda-2b^2\lambda+\lambda^2)\\
r_3&=-2(a^2-\lambda)(a^2b^2-c^2\lambda )^2\\
r_4&=(2c^2(a^2 b^2-(a^2+b^2)\lambda)^2\\
s_1&=-\frac{r_2}{r_1},\;\;\; s_2=-\frac{r_4}{r_3}
\end{align*}

Substituting $\cos\theta $ above for $g $ in \eqref{eqn:spatial-general} and obtain the spatial integral for the average cosine. 
Therefore it follows that
\begin{equation}\aligned
\int dx&= \int_0^{2\pi} \kappa_c^{\frac{2}{3}}ds  =\int_0^{2\pi} \frac{(a^2-\lambda)^{\frac{1}{3}}(b^2-\lambda)^{\frac{1}{3}}}{\sqrt{s_3} \sqrt{1-s_3\cos^2u}}du\\
   \int \kappa_c^{\frac{2}{3}}\cos\theta ds &= \int_0^{2\pi} \kappa_c^{\frac{2}{3}}\cos\theta du =s_4\int_0^{2\pi}  \frac{1-s_1\cos^2u}{(1-s_2\cos^2u)\sqrt{1-s_3\cos^2u}}du\\
    s_3&=\frac{a^2-b^2}{a^2-\lambda}, \;\;\;s_4=\frac{r_1}{r_3\sqrt{s_3}}(a^2-\lambda)^{\frac{1}{3}}(b^2-\lambda)^{\frac{1}{3}}
    \endaligned
\end{equation}

In terms of the elliptic integrals $K$ and $\Pi$ it follows that: 


\[ \overline{C}=\frac{r_1}{r_3}
\frac { \left(  {s_2} -s_1\right)  \Pi \left( s_2,\sqrt {s_3} \right) +{s_1}\,  K \left( 
\sqrt {s_3} \right)  }{ s_2 \, K \left(\sqrt {s_3}\right)  }
\]

In \cite{reznik2020-intelligencer,akopyan2020-invariants,bialy2020-invariants} the following expression was presented for the invariant sum of cosines in N-periodics:

\[ \sum{\cos\theta_i}=J L- N\]

\begin{figure}
    \centering
    \includegraphics[width=\textwidth]{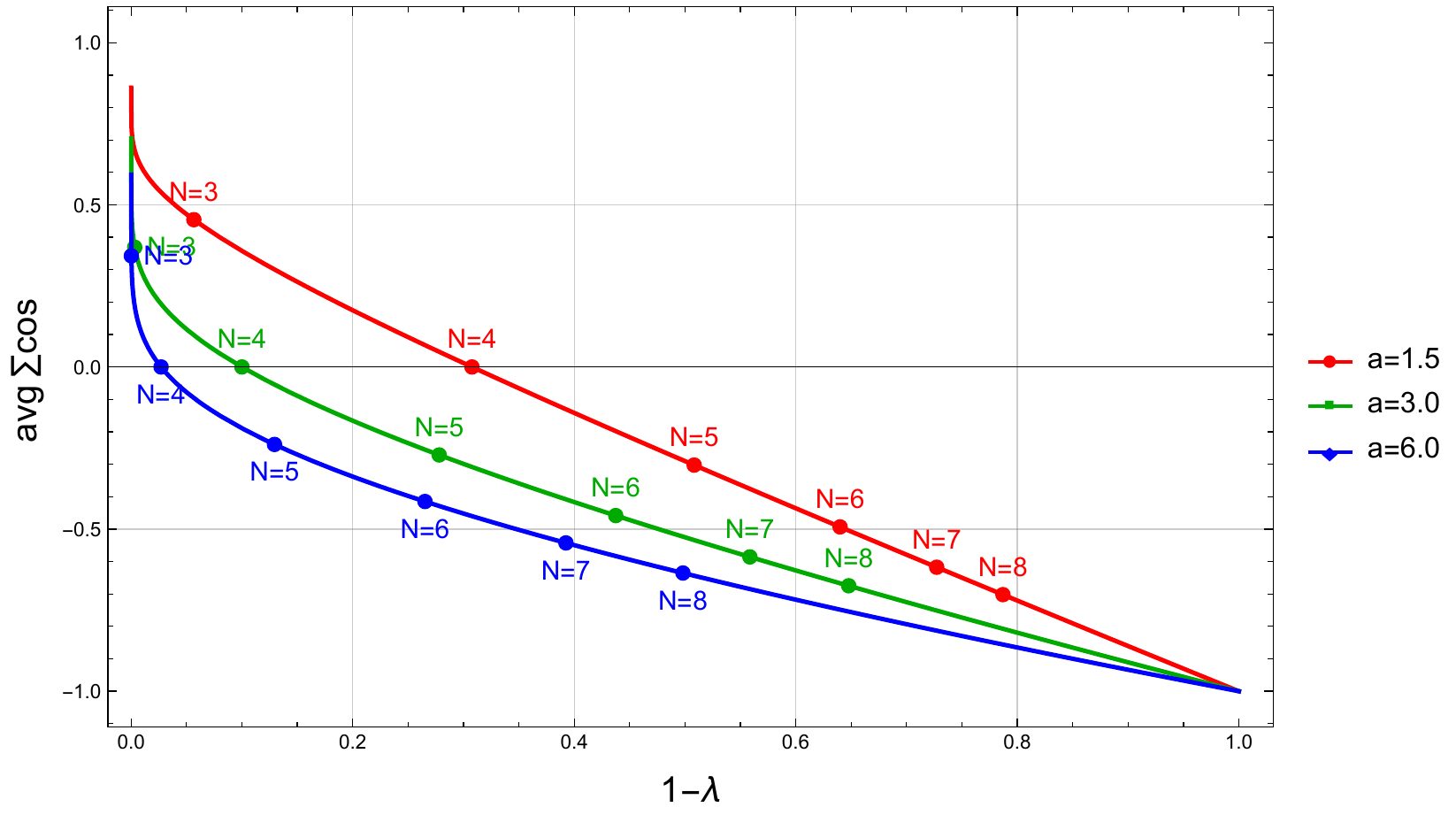}
    \caption{The average cosine vs $1-\lambda$, $b=1$ for three values of $a$, with $b=1$. The dots show agreement of the value with $J L/N-1$ for various non-intersecting N-periodics. When $1-\lambda$ is one (resp. zero), the average cosine tends to $1$ (resp. $-1$).}
    \label{fig:avg-cos-lambda}
\end{figure}

Therefore the average cosine for N-periodics is simply $JL/N-1$. \cref{fig:avg-cos-lambda} shows results obtained with spatial integration, and that they agree with the values predicted for N-periodics at the appropriate locations.

\subsubsection*{Sum of curvatures to two-thirds}

In \cite{reznik2021-invariants} we show conservation of $\sum{\kappa_i^{2/3}}$ is a corollary to the sum of cosines, where $\kappa_i$ denotes the curvature of the outer ellipse at the ith vertex. One can express $\kappa^{2/3}$ as a linear function of $\cos\theta$:

\[
\kappa^{2/3} =(a b)^{-\frac{4}{3}} \left(\frac{x^2}{a^4}+\frac{y^2}{b^4}\right)^{-1}=\frac{(ab)^{\frac{2}{3}}}{d_1 d_2}=\frac{4(a b)^{-\frac{2}{3}}}{|\nabla f|^2}=(ab)^{-\frac{4}{3}}\left( \frac{1+cos\theta}{2J^2}\right)-1
\]

\noindent Therefore, the sum of $\kappa^{2/3}$ is also invariant and its average value will be given by:

\[ \overline{\kappa^{2/3}}=   \frac{ 1}{\oint \,  \kappa_c^{2/3}  ds} \,  \, \oint \,   \kappa^{2/3}(s)\,  \kappa_c^{2/3} ds .
\]

%% file: 060_outer_cosines.tex
Referring to Figure~\ref{fig:definitions}, let $\theta_i'$ denote the ith internal angle of the outer polygon whose sides are tangent to the elliptic billiard at the vertices of an N-periodics. The product of $\theta_i'$ is invariant over N-periodics, for all N \cite{akopyan2020-invariants,bialy2020-invariants}. The geometric mean $\overline{C'}$ of $\theta_i'$ is given asymptotically by:

\[  \overline{C'} = \lim_{k \to \infty}\left(
 \prod_{i=1}^k   \cos\theta_i'\right)^{1/k} \]

To work with spatial integrals we must first convert the above to a sum:

\[ \log \overline{C'} =  \lim_{k \to \infty}(1/k)
 \sum_{i=1}^k   \log\left|\cos\theta_i'\right|
\]

As before, replace the above time average by the following spatial integral:

\[\log \overline{C'}=\frac{ 1}{\oint\,\kappa_c(s)^{2/3}ds}\,\,\oint \,\log|\cos\theta'(s)|\,\kappa_c^{2/3}ds\]

\begin{figure}
    \centering
\includegraphics[width=\textwidth]{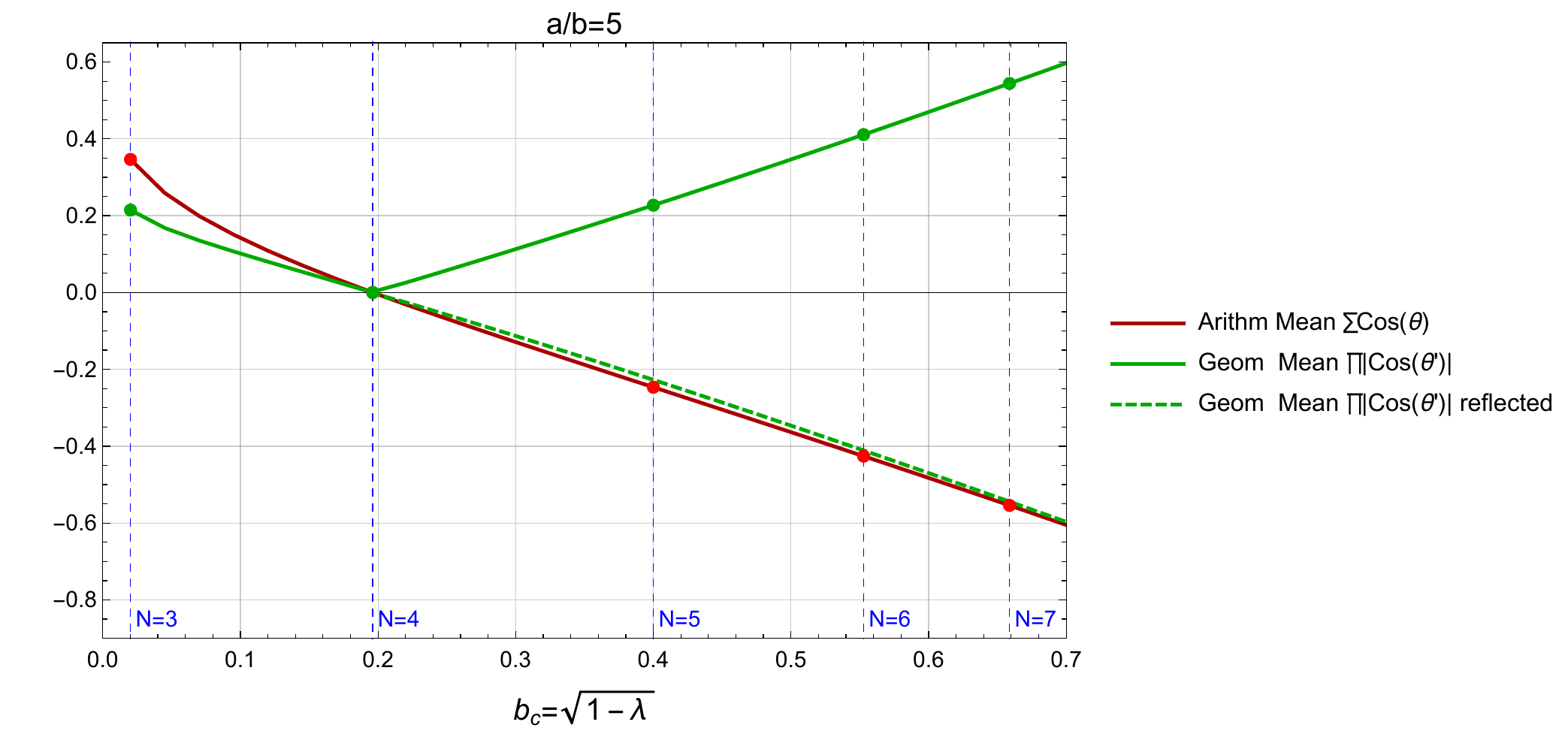}
    \caption{Average cosines (red) and geometric mean of outer cosines (green) vs. $b_c$, the minor semiaxes of the caustic. Here $a=5,b=1$. Blue dashed vertical lines mark the $b_c$ for non-intersecting orbits. Dashed green: past the $N=4$ caustic, the latter, the latter is reflected about the $x$ axis showing proximity to the average cosine.}
    \label{fig:koiller-sum-prod}
\end{figure}

A quick look  on  a picture shows that in order to compute $\cos \theta'$ it suffices to make the scalar product of the normalized gradients at the points $P_1,P_2$.

\begin{equation}
\cos \theta'  =   \frac{x_1 x_2/a^2 + y_1 y_2/b^2 }{  ( x_1^2/a^2 + y_1^2/b^2)^{1/2}   ( x_2^2/a^2 + y_2^2/b^2 )^{1/2}  }
\end{equation}

\begin{align*}
     \cos \theta'   &=- \frac{ c_a   \sqrt {-{a}^{2}{c}^{2} 
   \cos^2 u  + \left( {a}^{2}-\lambda \right)^{2}}}{\sqrt {{c}^{2}c_a^{2}   \cos^2 u  - \left( 2\,{b}^{2}\lambda+{
  c_a} \right) ^{2} \left( {a}^{2}-\lambda \right) }}
\\
 c_a &=a^2b^2-\lambda(a^2+b^2),\; \text{sign}(\cos\theta')=-\text{sign}(c_a)
\end{align*}

Numerical results for both average cosines and geometric mean of outer cosines are shown in \cref{fig:koiller-sum-prod} for $a=5$ (smaller $a$ make the two spatial averages become to close to each other). For values of $b_c$ where the trajectory is periodic, results obtained with spatial averages perfectly match numerically-estimated discrete averages computed numerically with N-periodics.

%% file: 080_questions.tex
The following questions are still unanswered:

\begin{itemize}
    \item Why is the geometric mean of outer aperiodic cosines so close to the average aperiodic cosines?
    \item Is there a universal measure expressed in terms of the outer ellipse?
    \item Can we use this framework to estimate aperiodic averages for cases where the caustic is a hyperbola?
    \item A third invariant introduced in \cite{reznik2020-intelligencer} was the ratio of outer-to-orbit areas. These do not seem amenable to a discrete sum of individual quantities. Would there be counterpart be for aperiodic areal averages?
\end{itemize}

%% file: 110_acknowledgments.tex
We would like to thank Sergei Tabachnikov and Arseniy Akopyan, and Hellmuth Stachel for invaluable discussions.

The second author is fellow of CNPq and coordinator of Project PRONEX/ CNPq/ FAPEG 2017 10 26 7000 508.